\documentclass{tran-l}
\usepackage{amssymb,latexsym, amscd}
\usepackage[all]{xy}
\usepackage{graphicx}

 \newtheorem{theorem}{Theorem}[subsection]
 
 \newtheorem{lemma}[theorem]{Lemma}
 
 \theoremstyle{definition}
 \newtheorem{definition}[theorem]{Definition}
 \theoremstyle{definition}
 
 \theoremstyle{remark}
 \newtheorem{rem}[theorem]{Remark}
 \numberwithin{equation}{subsection}

\usepackage{latexsym}
\usepackage{amssymb}
\newcommand{\ben}{\begin{equation}}
\newcommand{\een}{\end{equation}}

\newcommand{\integer}{\ensuremath{{\mathbb Z}}}
\newcommand{\naturals}{\ensuremath{{\mathbb N}}}
\newcommand{\real}{\ensuremath{{\mathbb R}}}
\newcommand{\complex}{\ensuremath{{\mathbb C}}}

\newcommand{\U}[1]{\ensuremath{{\mathrm U( #1 )}}}
\newcommand{\Or}[1]{\ensuremath{{\mathrm{O}( #1 )}}}

\newcommand{\Aa}{{\mathcal A}}

\newcommand{\UU}{{\mathcal U}}

\newcommand{\SSS}{{\mathcal S}}
\newcommand{\FF}{{\mathcal F}}

\newcommand{\LL}{\mathcal{L}}

\newcommand{\TT}{\mathcal{T}}

\newcommand{\Xx}{\mathsf{X}}
\newcommand{\Gg}{\mathsf{G}}

\newcommand{\target}{\mathsf{t}}
\newcommand{\source}{\mathsf{s}}

\newcommand{\Loop}{\mathsf{L}}

\newcommand{\To}{\longrightarrow}

\newcommand{\hol}{\mathrm{hol}}

\begin{document}

\title[Differential Characters over Orbifolds]
{Differential Characters on Orbifolds and String Connections I.}

\author{ Ernesto Lupercio and Bernardo Uribe}
\thanks{The first author was partially supported by the National Science
Foundation and Conacyt-M\'exico}

\address{Departamento de Matem\'{a}ticas, CINVESTAV, 
     Apartado Postal 14-740 
     07000 M\'{e}xico, D.F. M\'{E}XICO}
\address{Department of Mathematics, University of Michigan, Ann Arbor, MI 48109, USA} 
\email{ lupercio@math.wisc.edu \\ uribeb@umich.edu}

\begin{abstract} In this paper we introduce the Cheeger-Simons cohomology of a global quotient orbifold. We prove
that the Cheeger-Simons cohomology of the orbifold is isomorphic to its Beilinson-Deligne cohomology.
Furthermore we construct a string connection (\`{a} la Segal) from a global gerbe with connection
over the loop orbifold refining the corresponding differential character. \end{abstract}

\maketitle
\section{Introduction}

In our previous papers we have introduced both the concept of $n$-gerbe $L$ with
 connection over an orbifold $\Xx$, and the definition of orbifold Beilinson-Deligne
 cohomology. Our original motivation was the understanding of the concepts of $B$-field
 and discrete torsion in orbifold string theories \cite{LupercioUribe1,LupercioUribe3,
 LupercioUribe2} (cf. \cite{AdemRuan,Ruan2, AndoFrench, Ruan1,Sharpe}).  We showed that a $B$-field in mathematical
 terminology is the same as a 1-gerbe with connection over an orbifold (at the low
 energy limit, of course) and discrete torsion is a particular kind of flat $B$-field.
 String theory thus motivates the consideration of gerbes with connection over orbifolds.
 But the impact of considering gerbes goes beyond string theory. For example Y. Ruan \cite{Ruan3} 
 has recently implemented the usage of gerbes and string connections in order to obtain the twisted version of the Chen-Ruan quantum
 cohomology of an orbifold.

Given a smooth manifold $M$ and unitary line bundle with connection $(L,A)$ over $M$,
 we can consider consider its holonomy as a map
$$ \hol \colon Z_1(M) \To \U{1}.$$
We define  $\chi$ to be $$\chi:= -\frac{\sqrt{-1}}{2 \pi}\log \hol.$$ If we consider the curvature of $L$ as
 a 2-form $\omega$ on $M$ we have obtained a pair $(\chi, \omega)$ with
 $$\chi \colon Z_1 (M) \To \real / \integer$$
 and
 $$\chi(\partial c ) =  \int_c \omega \hbox{ \ mod \ } \integer$$
 whenever $c$ is a smooth 2-chain.

Following Cheeger-Simons \cite{CheegerSimons} we will denote by $\hat{H}^q _{cs}(M)$
 the $q^{th}$-group of differential characters of $M$.  The previous discussion with
 line bundles refers only to the case $q=2$. In the general case we have to substitute
 the line bundle by a $(q-2)$-gerbe with connection. The holonomy becomes now a homomorphism $Z_{q-1}(M) \to \real$.

Isomorphism classes of gerbes with connection are classified by the so-called
 Beilinson-Deligne cohomology \cite{Brylinski, LupercioUribe2}, denoted by
 $H^{q}(\Xx; \integer(q))$. It is a remarkable theorem the one that states that
$$H^{q}(\Xx; \integer(q)) \cong  \hat{H}^q _{cs}(M).$$
For this is stating that the holonomy and curvature of a gerbe completely determine its isomorphism class.

In section 2 of this paper we generalize the previous picture to the case of
 an orbifold of the form $\Xx=[M/G]$ (the previous manifold case being when $G=1$).
 One of the main difficulties is that of defining what we mean for a differential
 character on an orbifold. This is not completely obvious and there would be
 several possibilities for the definition. We offer one that allows a generalization
 of the main results to the orbifold setting. Our main construction actually works
 fine for orbifolds resulting for a \emph{Lie group} acting smoothly on a manifold with finite stabilizers.

In any case our main result is then
\begin{theorem}
For the orbifold $\Xx=[M/G]$,  the Beilinson-Deligne cohomology and the  Cheeger-Simons 
cohomology are canonically isomorphic.
\end{theorem}

In section 3 below we further refine this theorem. A gerbe with a connection
 contains more information that its isomorphism class, and therefore more information
 that just its differential character. The refinement we put forward in this paper
 is a \emph{string connection} associated to the gerbe with connection over the orbifold.
 Here we mean a string connection in the
 sense of Segal \cite{Segal02} and not in the closely related sense of Stolz-Teichner.
 In the manifold case ($G=1$) this can be encoded by  a functor from the category
 $\SSS^p(M)$, whose objects are maps from $(p-1)$-manifolds to $M$ and whose morphisms
 are maps of $p$-cobordisms to $M$, to the category of one-dimensional vector spaces. 

For example in the case of the line bundle the objects of $\SSS^1(M)$ are points in
 $M$ and its morphisms are paths between those points. The holonomy of a line bundle
 affords us a functor from $\SSS^1(M)$ to one-dimensional vector spaces. The vector
 space associated to a point is simply the fiber of the line bundle at that point. 

Going up one level a gerbe with connection produces via transgression a line
 bundle over the loop space of the manifold, and additionally a functor from
 $\SSS^2(M)$ to one-dimensional vector spaces by holonomy \cite{Murray02, LupercioUribe6}.
 In \cite{LupercioUribe6} we provided a generalization of the loop space for an orbifold
 called the \emph{loop groupoid} $\Loop \Xx$

The  second main result of this paper is the following.

\begin{theorem}
Let $\xi$ be a global gerbe with connection over $\Xx=[M/G]$ and $E$ the line
bundle with connection induced by it via transgression. Then $\xi$ permits to define
a string connection $\UU$ over the line bundle $E$ of the loop groupoid $\Loop \Xx$.  This string connection refines the corresponding differential character.
\end{theorem}

We would like to use this opportunity to thank the organizers of this excellent conference
 Profs. Tyler Jarvis,  Takashi Kimura   and   Arkady Vaintrob. We would like to thanks
 enlightening conversations with L. Borisov, T. Nevins, Y. Ruan and G. Segal. We dedicate
 this paper to Prof. Samuel Gitler on the occasion of his 70th birthday.

\section{Differential Characters}
In \cite{LupercioUribe6} we have put forward the definition of an $n$-gerbe with connective structure over 
an orbifold (i.e. \'{e}tale, proper foliation groupoid). We define them there as cocycles of a  cohomology
theory called  the Beilinson-Deligne (BD) cohomology. In this section we will define the orbifold
version of the Cheeger-Simons (CS) cohomology \cite{CheegerSimons} 
and then we will show that these two theories are isomorphic.

For smooth manifolds and certain algebraic varieties the isomorphism of this section has been proved by several authors 
\cite{BrylinskiMcLaughlin, DupontHainZucker, Esnault}. Here we are
concerned with the orbifold case and so, modifying Brylinski's definition \cite{Brylinski}
of BD-cohomology and Hopkins-Singer's \cite{HopkinsSinger} of CS-cohomology,
we will work out this theorem in the equivariant case, in which we have group actions.

\subsection{Cheeger-Simons Cohomology}

To make the exposition manageable we will start by considering orbifolds\footnote{We should not
 confuse the orbifold $\Xx=[M/G]$ with the quotient space $X=M/G$, we will remember all along 
the orbifold structure, including the stabilizers of  the $G$-action on $M$.} 
of the type $\Xx=[M/G]$ where $M$ is a finite dimensional, paracompact, smooth manifold
with a smooth action of a finite group $G$.  In this paragraph we will define the $G$-invariant
CS-cohomology and in the next section \ref{subsection equivariant CS} we will concentrate on the equivariant CS-cohomology.

\subsubsection{The case of a finite group $G$}

We want to define a cohomology theory of orbifolds $[M/G]$, which encodes the notion
of a $G$-invariant $q$-form with integral periods. If $M_G := M \times_G EG$ and $\Omega^q_{cl}(M)^G$
stands for the $G$-invariant closed $q$-forms over $M$, we want to
construct a theory that naturally fits in the upper-left corner of the diagram

 $$\xymatrix{
        ? \ar[r] \ar[d] & \Omega^q_{cl}(M)^G \ar[d] \\
        H^q(M_G; \integer) \ar[r] & H^q(M_G;\real).
        }$$

We will do this via the complex of smooth cochains. When $G=\{1\}$ this recovers the classical case
  due to Cheeger and Simons.
We will follow closely Hopkins and Singer's alternative definition \cite[Sect. 3.2]{HopkinsSinger}.

Let the complex $\hat{C}(q)^*(\Xx)$ be given by
\begin{eqnarray*}
\hat{C}(q)^n(\Xx) =
\left\{ 
\begin{array}{cc}
C^n(M\times EG;\integer)^G \times C^{n-1}(M \times EG ;\real)^G \times \Omega^q(M)^G & n \geq q \\
C^{n-1}(M \times EG;\integer)^G \times C^{n-2}(M\times EG ;\real)^G & n < q,
\end{array} \right.
\end{eqnarray*}
with differential 
$$d(c,h,\omega) := (\delta c, \tilde{\omega} - {c} - \delta h, d \omega)$$
$$d(c,h) = \left\{
\begin{array}{cc}
(\delta c, -{c} - \delta h, 0) & (c,h) \in \hat{C}(q)^{q-1}(\Xx) \\
(\delta c, -{c} -\delta h) & \mbox{otherwise}
\end{array} \right.
$$

An $n$-form $\omega \in \Omega^n(M)$  defines a smooth $n$-cochain in $C^n(M; \real)$ through integration; we will
use the same symbol for the form and the associated smooth cochain. Pulling the cochain back, via the projection map
$M\times EG \to M$, we get the $n$-cochain in $C^n(M \times EG ; \real)$ that we have called $\tilde{\omega}$; as
$G$ is finite $\tilde{\omega}$ defines a cochain in $C^{n}(M \times EG;\real)^G$.
By $C^{n}(M \times EG)^G$ we mean the $G$-invariant cochains\footnote{This construction cannot be generalized to the case when $G$ is a general Lie group. For when $G$ is a finite group $H^*(C^*(M\times EG;\integer)^G) \cong H^*(C^*(M\times_G EG; \integer))$, and we need the finiteness condition
in order to have a pushforward map at the cochain level $C^k(M \times EG)^G \to C^k(M\times_G EG)$.} .

\begin{definition}
The cohomology of the complex $\hat{C}(q)^n(\Xx)$ is the Cheeger-Simons cohomology of $\Xx$, i.e
$$\hat{H}(q)^*(\Xx) := H(\hat{C}(q)^*(\Xx))$$
\end{definition}

The homotopy Cartesian square 
$$\xymatrix{
        \hat{C}(q)^*(\Xx) \ar[r] \ar[d] & \Omega^{* \geq q}(M)^G \ar[d] \\
        C^q(M \times EG; \integer)^G \ar[r] & C^q(M \times EG;\real)^G
        }$$
yields a Mayer-Vietoris sequence
$$\cdots \to \hat{H}(q)^n(\Xx) \to H^n(M_G; \integer) \times H^n(\Omega^{* \geq q}(M)^G ) \to \ \ \ \  \ \ \ \ \ \ \ \ \ \ \ \ \ \ \ $$
$$ \ \ \ \ \ \ \ \ \ \ \ \ \ \ \ \ \ \ \ \ \ \ \ \ \ \ \ \ \ \ \ \ \ H^n(M_G ; \real )
\to  \hat{H}(q)^{n+1}(\Xx) \to \cdots$$
that leads to natural isomorphisms
$$\hat{H}(q)^n(\Xx) = \left\{ 
\begin{array}{cc}
H^n(M_G; \integer) & n >q \\
H^{n-1}(M_G; \real/\integer) & n<q,
\end{array} \right. $$
and a  short exact sequence
\begin{eqnarray} 
0 \to H^{q-1}(M_G, \real/\integer) \to \hat{H}(q)^q(\Xx) \stackrel{}{\to} \Omega^q_0(M)^G
 \to 0. \label{SES CS}
\end{eqnarray}
where $\Omega^j_0(M)^G$ stands for the $G$-invariant closed $j$-forms with integer periods.
This short exact sequence is obtained by taking $n=q$ in the previous long exact sequence.
We have that   $H^q(\Omega^{* \geq q}(M)^G )$
equals $\Omega^q_{cl}(M)^G$ and then one finds the quotient of the closed forms by the cokernel of the homomorphism
$H^q(M_G;\integer) \to H^q(M_G;\real)$ yielding the closed $q$-forms with integer periods.

\begin{rem}
 It is not entirely clear how to define differential characters analogous to the ones in   \cite{CheegerSimons} for the previous definition. 
If we want to do so we have to take a smooth representative of $M_G$ and define the differential forms on it.
This will allow us to generalize the concept of differential character to the equivariant setting. This will be done in the 
following section.
\end{rem}

\subsubsection{Equivariant Cheeger-Simons} \label{subsection equivariant CS}

Let $G$ be now a Lie group that acts smoothly on $M$ and let $\Or{n}$ be the orthogonal group on which $G$ embeds
as a subgroup. Let $V_{k+n,n}$ $(k\geq 0)$ be the Stiefel manifold of $n$-dimensional orthogonal frames in $\real^{k+n}$ on which
$\Or{n}$ acts freely (and \emph{a fortiori} $G$ also acts freely). Take $$M_G^k := M \times_G V_{k+n,n}$$
and let $M_G$ be the direct limit of the $M_G^k$'s. As $G$ acts freely in $M \times V_{k+n,n}$ then by the slice theorem the spaces 
$M_G^k$ are smooth manifolds and its limit is homotopy equivalent to the Borel construction $M\times_G EG$.

\begin{definition}
An equivariant differential $q$-form $\omega$ on $M_G$ is a sequence $\{\omega_k\}_{k \in \naturals}$ of differential
$q$-forms  $\omega_k \in \Omega^q(M^k_G)$ such that ${\rho_{k,k'}}^* \omega_{k'} =\omega_k$ with $k' \geq k$ and
$\rho_{k,k'}$ is the inclusion of manifolds $M^k_G \hookrightarrow M^{k'}_G$.
The vector space of such forms will be called $\Omega^q_G(M)$.

As the exterior derivative commutes with pullbacks, in our case $d ({\rho_{k,k'}}^* \omega_{k'})= {\rho_{k,k'}}^*(d  \omega_{k'}) $, 
then it makes sense to define the exterior derivative of $\omega$ as the sequence $d w = \{d \omega_k \}_{k \in \naturals}$.
Let 
$$\cdots\stackrel{d}{\to} \Omega^{q-1}_G(M) \stackrel{d}{\to} \Omega^{q}_G(M) \stackrel{d}{\to} \Omega^{q+1}_G(M) \stackrel{d}{\to} \cdots$$
be the equivariant De Rham complex of $M$. Its cohomology $H^*(\Omega^*_G(M))$ is the equivariant De Rham
cohomology of $M$, $H^*_{{\mathrm{DR}},G}(M)$.
It is well known \cite{Audin}
that there is a canonical isomorphism
$$H^*_{{\mathrm{DR}},G}(M) \cong H^*_G(M ; \real)$$
where $H^*_G(M ; \real) := H^*(M_G;\real)$ is the equivariant cohomology of $M$.
\end{definition}

We will write $C_{k}^G(M; R)$ ($C^k_G(M;\integer)$) for the $k$-chains (rep. $k$-cochains) on $M_G$ with
 coefficients in the ring $R$. Define the equivariant Cheeger-Simons complex $\hat{C}(q)^*_G(M)$ by

\begin{eqnarray*}
\hat{C}(q)_G^n(M) =
\left\{ 
\begin{array}{cc}
C^n_G(M;\integer) \times C^{n-1}_G(M;\real) \times \Omega^q_G(M) & n \geq q \\
C^{n-1}_G(M ;\integer) \times C^{n-2}_G(M;\real) & n < q,
\end{array} \right.
\end{eqnarray*}
with differential 
$$d(c,h,\omega) := (\delta c, \tilde{\omega} - {c} - \delta h, d \omega)$$
$$d(c,h) = \left\{
\begin{array}{cc}
(\delta c, -{c} - \delta h, 0) & (c,h) \in \hat{C}(q)_G^{q-1}(M) \\
(\delta c, -{c} -\delta h) & \mbox{otherwise}
\end{array} \right.
$$

An $n$-form $\omega \in \Omega_G^n(M)$ defines a smooth $n$-cochain $\tilde{\omega}$ in the following sense.
For $S \in C^G_n(M; \real)$ a $n$-chain there exist $k \in \naturals$ such that the image of $S$ in $M_G$ is included in the
subspace $M^k_G$ . Define $$\tilde{\omega} (S) := \int_S\omega_k.$$
It is clear from the definition that $\tilde{\omega}$ is independent of $k$.

\begin{definition}
The cohomology of the complex $\hat{C}(q)_G^n(M)$ is the equivariant Cheeger-Simons cohomology of $M$, i.e
$$\hat{H}(q)_G^*(M) := H(\hat{C}(q)_G^*(M))$$
\end{definition}

The homotopy Cartesian square 
$$\xymatrix{
        \hat{C}(q)_G^*(M) \ar[r] \ar[d] & \Omega^{* \geq q}_G(M) \ar[d] \\
        C^q_G(M ; \integer) \ar[r] & C^q_G(M ;\real)
        }$$
yields a Mayer-Vietoris sequence
$$\cdots \to \hat{H}(q)_G^n(M) \to H^n_G(M; \integer) \times H^n(\Omega_G^{* \geq q}(M) ) \to \ \ \ \  \ \ \ \ \ \ \ \ \ \ \ \ \ \ \ $$
$$ \ \ \ \ \ \ \ \ \ \ \ \ \ \ \ \ \ \ \ \ \ \ \ \ \ \ \ \ \ \ \ \ \ H^n_G(M ; \real )
\to  \hat{H}(q)_G^{n+1}(M) \to \cdots$$
that leads to natural isomorphisms
$$\hat{H}(q)_G^n(M) = \left\{ 
\begin{array}{cc}
H^n_G(M; \integer) & n >q \\
H^{n-1}_G(M; \real/\integer) & n<q,
\end{array} \right. $$
and a  short exact sequence
\begin{eqnarray} 
0 \to H^{q-1}_G(M, \real/\integer) \to \hat{H}(q)_G^q(M) \stackrel{\alpha}{\to} \Omega^q_{G,0}(M)
 \to 0. \label{SES E-CS}
\end{eqnarray}
where $\Omega^j_{G,0}(M)$ stands for the equivariant closed $j$-forms with integer periods, 
(a form $\omega$ has integer periods if $\tilde{\omega}(S) \in \integer$ whenever $\partial S = 0$).

Note that when the group $G$ is finite any $G$-invariant differential form
$\omega$ over $M$ pulls back to an equivariant differential form $\omega$ with $$\omega_k :=\frac{1}{|G|}\sum_{g \in G} g^*({\pi_k}^* \omega) $$ via the projection $\pi_k \colon M \times V_{k+n,n} \to M$.  
It is clear that the map $\gamma : \Omega^q(M)^G \to \Omega^q_G(M)$ is injective.
So if we consider the image of $ \gamma(\Omega^q_0(M)^G) $ in $\Omega^q_{0,G}(M)$ an then we consider the inverse image of this set
via the map $\alpha$ of the short exact sequence \ref{SES E-CS} we obtain the following short exact sequence: 
\begin{eqnarray} 
0 \to H^{q-1}_G(M, \real/\integer) \to \alpha^{-1}(\gamma(\Omega^q_0(M)^G)) \stackrel{\alpha}{\to} \Omega^q_{0}(M)^G
 \to 0. \label{SES CS invariant}
\end{eqnarray}
So, in view of \ref{SES CS} and \ref{SES CS invariant}  we get
\begin{lemma} \label{proposition G-invariant}
The CS-cohomology of the orbifold $[M/G]$ injects in the $G$-equivariant CS-cohomology of $M$. Moreover
$$  \hat{H}(q)^q([M/G]) \cong \alpha^{-1}(\gamma(\Omega^q_0(M)^G)).$$
Therefore the CS-cohomology of $[M/G]$ consist of the classes in the $G$-equivariant
CS-cohomology of $M$ that are obtained via a $G$-invariant closed differential forms over $M$ with integer periods. 
\end{lemma}

The previous description allows us to go one step further to define the {\bf equivariant differential characters} of $M$
generalizing the original construction of Cheeger and Simons \cite{CheegerSimons}
\begin{definition}
A $G$ equivariant differential character of $M$ of degree $q$ consists of a pair $(\chi, \omega)$ with 
$$\chi: Z_{q-1}(M_G ;\integer) \to \real / \integer$$
a character defined over the group of $q-1$-cycles, and $\omega \in \Omega^q_G(M)$ an equivariant differential $q$-form
such that for every smooth $q$-chain $S\in C_q(M_G;\real)$
$$\chi(\partial S) = \tilde{\omega}(S).$$

We will denote by $\hat{H}^q_{cs,G}(M)$ the $q^{th}$ equivariant group of CS differential characters of $M$.      
\end{definition}

As indicated in \cite{HopkinsSinger} the map 
\begin{eqnarray} \label{Iso CS G-difchar} 
\begin{array}{ccc} 
\hat{H}(q)^q_G(M) & \stackrel{\cong}{\to} & \hat{H}^q_{cs,G}(M) \\
(c,h,\omega) & \mapsto & (\chi, \omega)
\end{array}
\end{eqnarray}
where $\chi(z) := h(z) \hbox{  mod  } \integer$, is an isomorphism.

Following the spirit of  Proposition \ref{proposition G-invariant}
we can define the  {\bf differential characters of $[M/G]$} when $G$ is a finite group, as
the equivariant characters $(\chi, \omega)$ such that $\omega$ belongs to the image of some $G$-invariant
closed form with integer periods $\omega \in\Omega^q_0(M)^G$
 under the map $\gamma$.
\begin{definition}
The $q^{th}$ group of differential characters over the orbifold $[M/G]$  will
be denoted by $\hat{H}_{cs}^q([M/G])$.
\end{definition}

As before we get an isomorphism
\begin{eqnarray}\label{Iso CS difchar orbifold} 
\begin{array}{ccc}  
\hat{H}(q)^q([M/G]) & \stackrel{\cong}{\to} & \hat{H}^q_{cs}([M/G]) \\
(c,h,\omega) & \mapsto & (\chi, \gamma(\omega)).
\end{array}
\end{eqnarray}

\subsection{Beilinson-Deligne Cohomology}
BD-cohomology was discovered by Beilinson and Deligne for the purpose of having a cohomology theory for algebraic varieties
which includes singular cohomology and the intermediate Jacobians of Griffiths. We will deal with a  smooth analog of this
theory.

Recall that for a $\Xx$-sheaf , where $\Xx=[M/G]$, we mean a sheaf $\FF$ over $M$ on which $G$ acts continuously.
If $\FF$ is abelian, the cohomology groups $H^n(\Xx;\FF)$ are defined as the cohomology groups of the complex
$$\Gamma(M, \TT^0)^G \to \Gamma(M;\TT^1)^G \to \cdots$$
where $\FF \to \TT^0 \to \TT^1 \to \cdots$ is a resolution of $\FF$ by injective $\Xx$ sheaves and 
$\Gamma(M;\TT^j)^G$ are the $G$-invariant sections. When the abelian
sheaf $\FF$ is locally constant (for example $\FF = \integer$) is a result of Moerdijk \cite{Moerdijk98} that $H^*(\Xx; \FF)
\cong H^*(B\Xx; \FF)$ where the left hand side is sheaf cohomology and the right hand side is simplicial cohomology of
$B\Xx \simeq M_G$ with coefficients in $\FF$. 

Let $\Aa^p_{\Xx}$ denote the $\Xx$-sheaf of differential $p$-forms and $\integer_{\Xx}$ the constant $\integer$ valued $\Xx$
sheaf with $\integer_{\Xx} \to \Aa^0_{\Xx}$ the natural inclusion of constant into smooth functions.

\begin{definition}
The smooth BD complex $\integer(q)$ is the complex of $\Xx$ sheaves
$$\integer_{\Xx} \to \Aa^0_{\Xx} \stackrel{d}{\to} \Aa^1_{\Xx} \stackrel{d}{\to} 
\cdots  \stackrel{d}{\to}  \Aa^{q-1}_{\Xx},$$
and the hypercohomology groups $H^*(\Xx, \integer(q))$ are called the {\bf smooth Beilinson-Deligne cohomology}
of $\Xx$.
\end{definition} 

Now, let $\U{1}(q)$ be the complex of sheaves
$$\U{1}_{\Xx} \stackrel{\sqrt{-1}d\log}{\To} \Aa^1_{\Xx} \stackrel{d}{\to} 
\cdots  \stackrel{d}{\to}  \Aa^{q-1}_{\Xx}$$
where $\U{1}_{\Xx}$ is the sheaf of $\U{1}$-valued functions.
Because of the quasi-isomorphism between $\integer(q)$ and $\U{1}(q)[-1]$, i.e.
\begin{eqnarray}
\label{quasiisomorphism}
\xymatrix{
    \integer(p)_\Xx \ar[r] &
     \Aa^0_{\Xx} \ar[r]^{d} \ar[d]_{\exp(-i\_)} &\Aa^1_{\Xx} \ar[r]^{d} \ar[d] &
\cdots \ar[r]^{d} & \Aa^{q-1}_{\Xx} \ar[d]\\
    &    \U{1}_\Xx \ar[r]^{\sqrt{-1}d \log} &\Aa^1_{\Xx} \ar[r]^{d} &
\cdots \ar[r]^{d} &\Aa^{p-1}_{\Xx}
    }
\end{eqnarray}
there is an isomorphism of hypercohomologies
\begin{eqnarray} \label{isodelignecohomology}
H^{n-1}(\Xx, \U{1}(q)) \cong H^n(\Xx, \integer(p)).
\end{eqnarray}

We need to use a more computational approach to this cohomology theory, basically because we will be using
3-cocycles in order to define a string connection, and so we will use a \v{C}ech description of the
BD-cohomology. In order to make the exposition less lengthy, we are going to make use
of some results that can be found in our previous paper \cite{LupercioUribe6}.
As $M$ is paracompact, for the orbifold $\Xx =[M/G]$ ( or better, the proper \'{e}tale foliation groupoid
with objects $\Xx_0 =M$ and morphisms $\Xx_1 = M \times G$) we can find a smooth \'{e}tale 
Leray groupoid $\Gg$ together with
a Morita map $\Gg \to \Xx$, making $\Gg$ and $\Xx$ Morita equivalent. Being Leray means that the spaces $\Gg_n$ of
$n$-composable morphisms of $\Gg$ are diffeomorphic to a disjoint union of contractible open sets.
In the case when $G=\{1\}$ (i.e. $\Xx = M$)
this amounts to finding a contractible open cover of $M$ such that all the finite intersections
of this cover are either contractible or empty and then making $\Gg_n$ to be the disjoint union
of all intersections of $n$ sets in the cover.

Let's denote by
$\breve{C}^*(\Gg;\U{1}(q))$ the total complex
\begin{eqnarray*}
\xymatrix{
\breve{C}^0(\Gg; \U{1}(q))  \ar[r]^{\delta -d } &
         \breve{C}^1(\Gg;\U{1}(q) )
        \ar[r]^{\delta + d}  &  \breve{C}^2(\Gg; \U{1}(q))
         \ar[r]^{\ \ \delta-d} & \cdots  
}
\end{eqnarray*}
induced
by the double complex
\begin{eqnarray} \label{double complex}
   \xymatrix{
     \vdots & \vdots & \vdots & & \vdots\\
     \Gamma(\Gg_2, \U{1}_\Gg) \ar[u]^\delta \ar[r]^{\sqrt{-1}d \log} &
         \Gamma(\Gg_2, \Aa^1_{\Gg})
        \ar[r]^d \ar[u]^\delta &  \Gamma(\Gg_2, \Aa^2_{\Gg})
        \ar[u]^\delta \ar[r]^d & \cdots \ar[r]^d & \Gamma(\Gg_2, \Aa^{q-1}_{\Gg})
         \ar[u]^\delta\\
      \Gamma(\Gg_1, \U{1}_\Gg) \ar[u]^\delta \ar[r]^{\sqrt{-1}d \log} &
          \Gamma(\Gg_1, \Aa^1_{\Gg})
        \ar[r]^d \ar[u]^\delta &  \Gamma(\Gg_1, \Aa^2_{\Gg})
        \ar[u]^\delta \ar[r]^d & \cdots \ar[r]^d & \Gamma(\Gg_1, \Aa^{q-1}_{\Gg})
        \ar[u]^\delta\\
      \Gamma(\Gg_0, \U{1}_\Gg) \ar[u]^\delta \ar[r]^{\sqrt{-1}d \log} &
         \Gamma(\Gg_0, \Aa^1_{\Gg})
          \ar[r]^d \ar[u]^\delta &  \Gamma(\Gg_0, \Aa^2_{\Gg})
          \ar[u]^\delta \ar[r]^d & \cdots\ar[r]^d & \Gamma(\Gg_0, \Aa^{q-1}_{\Gg})
          \ar[u]^\delta
    }
\end{eqnarray}
with $(\delta +(-1)^{i} d )$ as coboundary operator, where the $\delta$'s are the maps induced
simplicial structure of the nerve of the category $\Gg$ and $\Gamma(\Gg_i, \Aa^j_{\Gg})$ stands
for the global sections of the sheaf that induces $\Aa^j_{\Gg}$ over $\Gg_i$ (see \cite{LupercioUribe6}).
Then the {\v{C}}ech hypercohomology of the complex of sheaves
$\U{1}(q)$ is defined as the cohomology of the
{\v{C}}ech complex $\breve{C}(\Gg;\U{1}(q))$:
$$\breve{H}^*(\Gg; \U{1}(q)):=H^*\breve{C}(\Gg;\U{1}(q)).$$

As the $\Gg_i$'s are diffeomorphic to a disjoint union of
contractible sets -- Leray -- then the previous cohomology
actually matches the hypercohomology of the complex
$\U{1}(q)$, so we get
\begin{lemma} \label{cech=hyper}
The cohomology of the \u{C}ech complex
$\breve{C}^*(\Gg,\U{1}(q))$ is isomorphic to the
hypercohomology of the complex of sheaves $\U{1}(q)$ and as $\Gg \to \Xx$ are isomorphic, then
$$\breve{H}^*(\Gg, \U{1}(q)) \stackrel{\cong}{\to}
H^*(\Gg;\U{1}(q)) \cong H^*(\Xx;\U{1}(q)) . $$
\end{lemma}

As we are only interested in the case $\Xx =[M/G]$ we can make a 
more explicit description of the Leray groupoid $\Gg$. Take a contractible open cover $\{U_i\}_{i \in I}$
of $M$ such that all the finite intersections of the cover are either contractible or empty, and with
the property that for any $g \in G$ and any $i \in I$ there exists $j \in I$ so that $U_ig=U_j$. Define
$\Gg_0$ as the disjoint union of the $U_i$'s with $\Gg_0 \stackrel{\rho}{\to} M=\Xx_0$ the natural map.
Take $\Gg_1$ as the pullback square
  $$\xymatrix{ \Gg_1 \ar[r] \ar[d] & M \times G \ar[d]^{\source \times \target} \\
             \Gg_0 \times \Gg_0 \ar[r]^{\rho \times \rho}& M \times M}$$
where $\source(m,g) =m$ and $\target(m,g) = mg$. This defines the proper \'{e}tale Leray groupoid $\Gg$ and
by definition it is Morita equivalent to $\Xx$.

\begin{lemma}
There is a natural short exact sequence
\begin{eqnarray}
0 \to \breve{H}^{q-1}(\Gg;\real/\integer) \stackrel{\sigma}{\to} \breve{H}^{q-1}(\Gg; \U{1}(q)) 
\stackrel{\kappa}{\to} \Omega^q_0(M)^G \to 0. \label{SES BD}
\end{eqnarray}
\end{lemma}
\begin{proof}
The map $\sigma$ is obtained by the inclusion of the locally constant $\real/\integer$-valued $\Gg$-sheaf into $\U{1}_{\Gg}$, 
it follows that $j$ is  injective.
 Now let's consider an element $[f] \in \breve{H}^{q-1}(\Gg; \U{1}(q))$.
It will consist of the ${q}$-tuple $(\theta_0, \dots, \theta_{q-1})$ with $\theta_0 \in \Gamma(\Gg_{q-1},\U{1}_\Gg)$
and $\theta_i \in \Gamma(\Gg_{q-1-i}, \Aa^i_\Gg)$ that satisfies the cocycle condition 
$d\theta_i + (-1)^{q-1}\delta \theta_{i+1} =0$.

From the construction of $\Gg$ we see that we can think of $\Gg_1$ as the disjoint union of all the intersections
of two sets on the base times the group $G$, i.e. 
$$\Gg_1 = \left( \bigsqcup_{(i,j) \in I\times I} U_i \cap U_j \right) \times G$$
where the arrows in $U_i \cap U_j \times\{g\}$ start in $U_i|_{U_j}$ and end in $(U_j|_{U_i})g$.

By the cocycle condition we know that
$$g^* \theta_{q-1}|_{U_{j}g} - \theta_{q-1}|_{U_{i}} = d \theta_{q-2}|_{U_{ij}\times\{g\}} \mbox{ \ \  in \ \ } U_{ij}$$
where $U_{ij} = U_i \cap U_j$.  So if we define $q$-forms $\omega_i$ locally by $\omega_i := d \theta|_{U_i}$ is easy
to see that once all are glued together they will induce a global $q$-form $\omega$ 
over $M$ which is $G$-invariant. The globality is obtained by
taking $g=1$ and noting that $\omega_i$ and $\omega_j$ agree in the intersection and the invariance is easily seen by taking $i=j$.
As $\omega$ is defined locally by exact forms then it follows that $\omega$ is exact.
We define $\kappa([f]) :=\omega$; it is well defined because if $f'= (\theta'_0, \dots, \theta'_{q-1})$ is cohomologous
to $f$ then $\theta'_{q-1} - \theta_{q-1}$ is exact, so $f$ and $f'$ define the same $q$-form.

We are now left to prove that $\kappa$ is surjective and that $\ker( \kappa) \subset {\mathrm{Im}}(\sigma)$. We will do
so by looking at the double complexes used in  the proof of the De Rham theorem and at the one
by the \v{C}ech description of the complex of sheaves $\integer(q)$. 
Recall that if  $\real_\Gg$ is the $\Gg$-sheaf of locally constant $\real$-valued functions then we know
that the complex \cite{Weil52}
$$\Aa^0_\Gg \stackrel{d}{\to} \Aa^1_\Gg  \stackrel{d}{\to} \dots $$
is a resolution of injective sheaves.

If we have a BD class $[\theta_0, \dots , \theta_{q-1}]$ as before, such that its image under $\kappa$ is zero, i.e $\omega=0$,
then the $q-1$-form given by $\phi_{q-1}$ is closed. As the groupoid is Leray, by a successive application
of the Poincar\'{e} lemma, we can find a chain $(\alpha_0 ,\dots, \alpha_{q-2} ) \in \breve{C}^{q-2}(\Gg;\U{1}(q))$
such that 
$$ (\theta_0, \dots , \phi_{q-1}) + (d + (-1)^{q-2}\delta)(\alpha_0 ,\dots, \alpha_{q-2} ) = (\theta'_0, 0, \dots, 0).$$
Then  $\theta'_0$ is locally constant  (because $d \log \theta'_0 =0$) and $\delta \theta'_0 =1$, so it
defines a \v{C}ech cocycle with values in the $\real/\integer$ $\Gg$ sheaf. This implies that the kernel of $\kappa$
is included in the image of $\sigma$.

Now, A  $G$-invariant $q$-form with integer periods $\omega$, via the De Rham theorem,  defines 
forms $\phi_i \in \Gamma(\Gg_{q-1-i}, \Aa^i_\Gg)$ and
a cocycle 
in $c \in\Gamma(\Gg_q, \real_\Gg)$ such that $d\phi_i + (-1)^{q-1}\delta \phi_{i+1} =0$, $\delta \phi_0 +(-1)^{q-1}c =0$
 and $\delta c =0$ (here we are making use of the quasi-isomorphism of \ref{quasiisomorphism}).
As $\omega$ has integer periods then there exist $c' \in \Gamma(\Gg_q, \integer_\Gg)$ and $h \in \Gamma(\Gg_{q-1}, \real_Gg)$
such that $c' = \delta h  + c$, then $(c, \phi_0 + (-1)^{q-2}h, \phi_1, \dots, \phi_{q-1})$  is a 
BD cocycle for the complex of $\Gg$
sheaves $\integer(q)$. Its BD-cohomology class under the map $\kappa$ is $\omega$. So $\kappa$ is surjective.

The sequence is short exact.

\end{proof}
 
As the groupoids $\Gg$ and $\Xx$ are Morita equivalent\footnote{Two such groupoids are Morita equivalent if and only if they represent the same orbifold (see for example \cite{Moerdijk2002}).} then we have the short exact sequence
\begin{eqnarray} \label{SES BD2}
0 \to H^{q-1}(\Xx;\real/\integer) \to H^{q}(\Xx; \integer(q)) \to \Omega^q_0(M)^G \to 0.
\end{eqnarray}

\begin{theorem}
For the orbifold $\Xx=[M/G]$, with $G$ a finite group, the Beilinson-Deligne cohomology and the  Cheeger-Simons 
cohomology are canonically isomorphic.
\end{theorem}
\begin{proof}
In view of the short exact sequences \ref{SES CS} and \ref{SES BD2} is just a matter of constructing
a map from $H^{q}([M/G]; \integer(q)) $ to $\hat{H}(q)^q([M/G])$. This turns out to be somewhat subtle in the case
of orbifolds and it is actually given by  the string connection described in the next chapter. For now we can bypass this by a careful  use of our definition
of equivariant CS-cohomology. 

What we will actually do is to define a map from $H^{q}([M/G]; \integer(q)) $ 
to the group of differential characters of $[M/G]$, namely
$\hat{H}^q_{cs}([M/G])$ (see \ref{Iso CS difchar orbifold} ). It consists of pairs $(\chi, \gamma(\omega))$
with $\omega \in \Omega^q_0(M)^G$ and 
$$\chi \colon Z_{q-1}(M_G ; \integer ) \to \real/\integer$$
such that for any smooth $q$ chain $z$ we have
$$\chi(\partial z) = \widetilde{\gamma(\omega)}(z) \hbox{ \ mod \ } \integer.$$

For an element $[\xi] \in H^q([M/G] ; \integer(q))$  via \ref{SES BD2} we obtain a form $\omega \in \Omega^q_0(M)^G$.
So we only need to define $\chi$ and we will do so by defining its value  on the $(q-1)$-dimensional, compact boundaryless submanifolds
of $M_G$. 
Then let $\Sigma$ be a compact $(q-1)$-dimensional smooth manifold without boundary and $\phi  \colon \Sigma \to M_G$ a smooth map.
Then there exist $k \in \naturals$ such that $\phi  \colon \Sigma \to M^k_G$ . As the group $G$ is finite
we can pullback the $G$-bundle $\pi_k \colon M \times V_{k+n,n} \to M \times_G V_{k+n,n}$ via $\phi$ and we call such $G$-bundle
 $P$, i.e. $P$ is a $G$-bundle over $\Sigma$ and is the pullback square of the following diagram
  $$\xymatrix{
     P \ar[rr] \ar[d] & &M \times V_{k+n,n} \ar[d]^{\pi_k} \\
     \Sigma \ar[rr]^\phi & & M \times_G V_{k+n,n}.
   }$$

Composing $\phi$ with the projection map $pr \colon M \times V_{k+n,n} \to M$ we obtain a $G$-equivariant
map $\hat{\phi} := pr \circ \phi \colon P \to M$. Note that this map can also be seen as a map 
of orbifolds $\hat{\phi} \colon [P/G] \to [M/G]$).
Pulling back $[\xi]$ via $\hat{\phi}$ we obtain a class $[\hat{\phi}^*\xi]$ in $H^q([P/G]; \integer(q))$.
As $\Omega^q(P)^G = \Omega^q(\Sigma)=\{0\}$ then the class $[\hat{\phi}^*\xi]$ is isomorphic to a
class $[\rho] \in H^{q-1}([P/G] ; \real/\integer) = H^{q-1}(\Sigma ; \real/\integer)$.
Define
$$\chi(\Sigma) : = \rho(\Sigma).$$
It is well defined because if $\rho' = \rho + \delta \kappa$ with $\kappa \in C^{q-2}(\Sigma ; \real )$, as $\partial \Sigma = 0$
then $\rho'(\Sigma ) = \rho(\Sigma) + \delta ( \kappa) (\Sigma)= \rho(\Sigma) + \kappa(\partial \Sigma) = \rho (\Sigma) $.

Taking $[\xi] \mapsto (\chi, \gamma(\omega))$ we get a map 
$$H^{q}([M/G]; \integer(q)) \to \hat{H}^q_{cs}([M/G])$$
that commutes with the short exact sequences \ref{SES CS} and \ref{SES BD2}. By the 5-lemma
the result follows.\footnote{ In the case when $G=\{1\}$ and the orbifold is simply a smooth manifold,
 the theorem was previously obtained in \cite{Murray02}.}
\end{proof}

\section{String Connections}\label{section string connection}

In this section we will focus our attention to gerbes with connection over the orbifold $[M/G]$
although most of the constructions could be generalized to $n$-gerbes with connection.

In order to make the exposition clearer let us introduce this section by 
explaining the idea of a string connection in the case of a manifold $M$ without any group action.

A gerbe with connection over $M$ induces a line bundle $E$ with connection over the free loop space $\LL M$ of
$M$ via a transgression map. The connection over $\LL M$ allows one to do parallel transport over a path 
$\gamma \colon [0,1] \to \LL M$  defining an invertible linear operator $A_\gamma$ between
the fibers $E_{\gamma(0)} $ and $E_{\gamma(1)}$. This path $\gamma$ in the loop space could be seen
also as a 2 dimensional submanifold of $M$ of genus zero  with boundary components 
the loops $\gamma(0)$ and $\gamma(1)$.
But we would like to do more than just being able to do parallel transport over a tube, we would like to
do a more general transport through an embedded oriented Riemann surface with boundary in $M$. This will give
us a transport  operator 
$$\UU_\Sigma: E_{\gamma_1} \otimes  \cdots \otimes E_{\gamma_p} \to 
E_{\gamma_{p+1}} \otimes  \cdots \otimes E_{\gamma_{p+q}} $$
 to each smooth surface  $\Sigma$ in $ M$ which has $p$ incoming parametrized boundary circles
$\gamma_1, \dots, \gamma_p$  and $q$
outgoing parametrized outgoing circles $\gamma_{p+1}, \dots, \gamma_{p+q}$. The case $p=q=1$  and $\Sigma$ a torus
is depicted below.
\begin{eqnarray}
\includegraphics[height=1.0in]{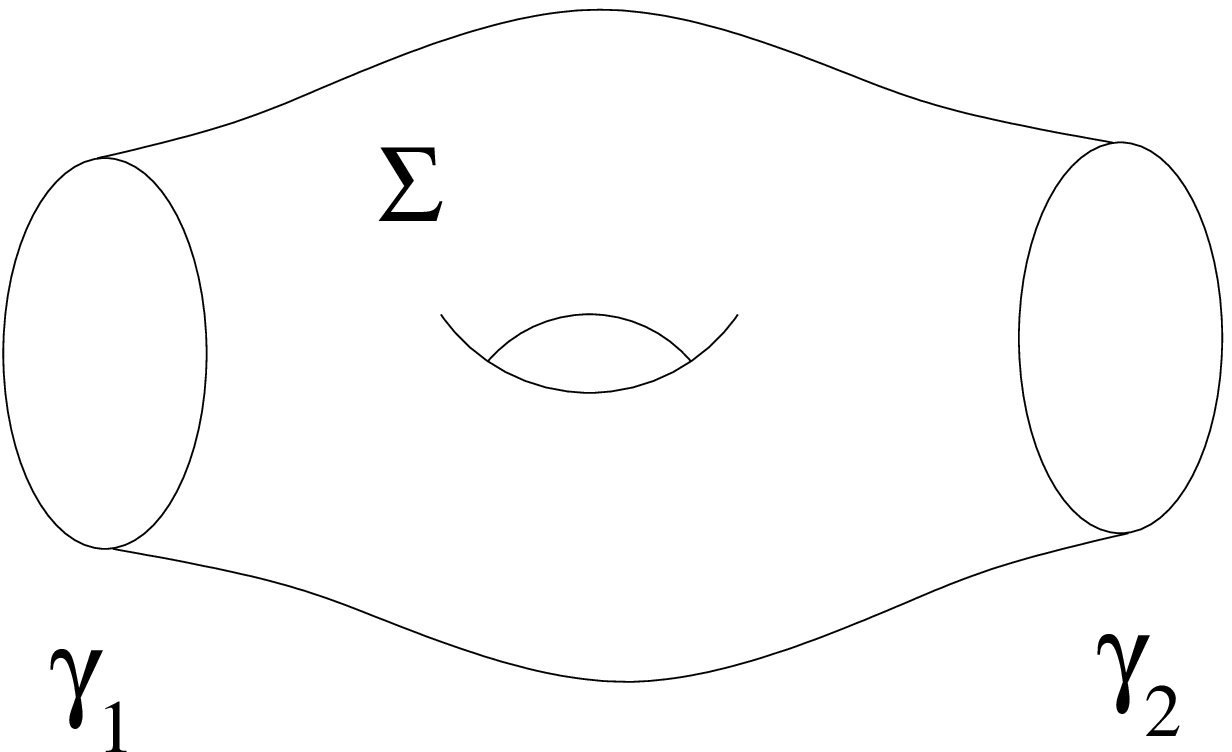} \label{graph string}
\end{eqnarray}
As for an ordinary connection $\gamma \mapsto A_\gamma$, the properties of the assignment $\Sigma \mapsto \UU_\Sigma$
are that is transitive with respect to concatenating Riemann surfaces, and that it is parametrization independent
in the sense that it does not change if $\Sigma \to M$ is replaced by the composite $\Sigma' \to \Sigma \to M$
where $\Sigma' \to  \Sigma$ is a diffeomorphism. Actually this is also true for self-gluing of an incoming and an outgoing components in the same Riemann surface.

A gerbe with connection over $M$ allows one to define such operator $\UU$. In this case the closed
3-form $\omega$ with integer periods (the field strength, or curvature) must be such that 
the value of $\omega(v)$ on an element $v$ of 3-volume at $x$ bounded by a surface $\Sigma$ is
given by 
\begin{eqnarray} \label{Segal's formula}
\exp(  \sqrt{-1} \ \omega(v)) = \UU_\Sigma,
\end{eqnarray}
where $\Sigma$ is regarded as a path in the loop space from a point loop at $x$ to itself (therefore
$\UU_\Sigma$ is an invertible $\complex$-linear map from the fiber at the constant loop $x$ to itself, hence
we can associate it a complex number of norm one).
 This operator $\UU$ is what Segal has called a ``String Connection''
\cite{Segal02}.

If we were to define $$\UU_\Sigma : = \exp(\sqrt{-1} \chi (\Sigma))$$
when $\partial v=\Sigma$ 
where $(\chi, \omega)$ is the differential character given by the gerbe with connection
we can see that formula \ref{Segal's formula} follows from the fact that
$$\chi(\Sigma) = \omega(\Sigma) \  \hbox{mod} \ \integer.$$
 So the gerbe with connection allows
one to define the operator $\UU$ on surfaces that are boundaries, but in order to do it for any
other surface we need to do more.

Coming back to orbifolds,  we will explain how the operator $\UU$ is defined for an specific 
gerbe with connection over $[M/G]$; we will call this class of gerbes {\bf global}.
\begin{definition}
 A global gerbe with connection $\xi$ 
is a Beilinson-Deligne cocycle over $[M/G]$ whose data is given by global forms\footnote{General gerbes require Leray representatives that are Morita equivalent to the given group action \cite{LupercioUribe6}.} . Namely $\xi$ will consist
of the forms  $B \in \Omega^2(M)$, $A_g \in \Omega^1(M)$ and $\rho_{g,h} \colon M \to \U{1}$ for
$g,h \in G$ such that it is a cocycle in the double complex \ref{double complex}, i.e.
\begin{eqnarray}
g^* B - B &= &d A_g \label{relation B}\\
A_g + g^* A_h - A_{gh} & = & \sqrt{-1} d \log \rho_{g,h}. \label{relation A} 
\end{eqnarray}
The curvature $\omega \in \Omega^3_0(M)$ of $\xi$ is $G$-invariant and we know that $\omega = dB$.
\end{definition}
A gerbe with connection $\xi$ induces a complex line bundle over the loop orbifold of $[M/G]$ and
the transport operator will act on its fibers. Here
we need to recall the definition of the loop orbifold, and as the gerbe we have in mind
is global it suffices to take orbifold
maps of principal bundles over the circle to $[M/G]$. The general definition of the loop orbifold
can be found in \cite{LupercioUribe2} and the definition of the line bundle induced by $\xi$ is in
\cite{LupercioUribe6}.

A loop over the orbifold $[M/G]$ will consist of a map $\phi \colon Q \to M$ of a $\Gamma$-principal
bundle $Q$ over the circle $S^1$ and a homomorphism $\phi_\# \colon \Gamma \to G$ such that
$\phi$ is $\phi_\#$-equivariant. Let's call this space of loops by $\LL [M/G]$. It has a natural 
action of the group $G$ as follows. For $h \in G$ let $\psi := \phi \cdot h$ where
$\psi(x) := \phi(x)h$ and $\psi_\#(\tau) = h^{-1} \phi_\#(\tau) h$, then $\psi \colon Q \to M$ and
is $\psi_\#$ equivariant.

\begin{definition}
The groupoid given by $\left(\LL [M/G]\right)/G$ is what we call the loop orbifold.
\end{definition}

Note that we could have just taken only $\Gamma$ principal bundles where $\Gamma$ is a finite
cyclic group $\integer_m$. This because the relevant information of such maps lie
on the holonomy of the circle on $[M/G]$, and this is characterized by a conjugacy class in $G$.
We did not do so because it will be needed in such generality to simplify what follows.

The topology of the loop groupoid is given by the compact-open topology of the space of maps
of a given principal bundle. Then  loops defined  over  two different principal bundles are in different
connected components of the loop orbifold.

The line bundle $E$ associated to the gerbe $\xi$ is obtained by defining a groupoid
map from  $\left(\LL [M/G]\right)/G$ to $\U{1}$. So it is a trivial line bundle $E = \complex \times \LL [M/G]$
over $\LL [M/G]$ with a $G$ action.
Define
\begin{eqnarray}
F \colon \LL [M/G] \times G & \to & \U{1} \label{def. F}\\
F(\phi,g)&\mapsto& \exp\left( \frac{\sqrt{-1}}{|\Gamma|} \int_Q \phi^* A_g\right) 
\end{eqnarray}
where $\phi \colon Q \to M$ is a loop with $Q$ a $\Gamma$-principal bundle. 
Using  equation  \ref{relation A} and the fact that  $\partial Q = 0$ it follows that
$$F(\phi,g)  F(\phi\cdot g, h) = F(\phi, gh),$$
meaning that the map $F$ is a map of groupoids, therefore defining an action on the  line bundle $E$.

Now we want to define a string connection over this line bundle $E$ given by the gerbe with connection $\xi$.
We need to consider the  equivalent over an orbifold of a Riemann surface with boundary.
This will consist of a map $\Phi : P \to M$ of a 
$\Gamma$-principal bundle $P$ over an oriented Riemann surface $\Sigma$ ($\Gamma$ finite) and a homomorphism
$\Phi_\# \colon \Gamma \to G$ such that $\Phi$ is $\Phi_\#$-equivariant. Note that there is a natural
action of the group $G$ on $\Phi$. It is defined in the same way as for loops.

Let the boundary of $\Sigma$ have $p$ incoming parametrized circles and $q$ outgoing. Then 
the boundary $\partial P$ of $P$ will consist of $p$ incoming orbifold loops $\gamma_i \colon Q_i \to M$ $1\leq i \leq p$
with the induced orientation,
and $q$ outgoing ones $\gamma_j \colon \overline{Q}_j \to M$, $p+1 \leq j \leq p+q$
with the opposite orientation so that $\partial P = \bigsqcup_i Q_i \sqcup \bigsqcup_j Q_j$.
Here the  $Q_i$'s and the $Q_j$'s are $\Gamma$-principal bundles
over the circle. As our line bundle is trivial (without the $G$ action) the operator
$\UU_{\Phi}$ is just a complex number.
Define
\begin{eqnarray*}
\UU_\Phi  \colon  E_{\gamma_1} \otimes \cdots \otimes E_{\gamma_p} = \complex
 & \to & \complex =E_{\gamma_{p+1}} \otimes \cdots \otimes E_{\gamma_{p+q}} \\
\UU_\Phi &:=& \exp\left( \frac{\sqrt{-1}}{|\Gamma|} \int_P \Phi^* B \right). 
\end{eqnarray*}

The only thing left to prove for $\UU$  to be a string connection over $E$ is that it
is compatible with the $G$ action. The concatenation property and the invariance
under diffeomorphisms that fix the boundary follow
from the fact that the operator is defined through an integral over the surface $P$.
 So we want the following diagram to be commutative for any $g \in G$
$$
  \xymatrix{
      \bigotimes_{i=1}^p E_{\gamma_i} \ar[rrr]^{\exp\left( \frac{\sqrt{-1}}{|\Gamma|} \int_P \Phi^* B \right)} 
      \ar[dd]_{\prod_i\exp\left( \frac{\sqrt{-1}}{|\Gamma|} \int_{Q_i} \gamma_i^* A_g\right)} 
      & & & \bigotimes_{j=p+1}^{p+q} E_{\gamma_j}
    \ar[dd]^{\prod_j\exp\left( \frac{\sqrt{-1}}{|\Gamma|} \int_{\overline{Q}_j} \gamma_j^* A_g\right)}   \\
  & &  & \\
      \bigotimes_{i=1}^p E_{\gamma_ig} \ar[rrr]^{\exp\left( \frac{\sqrt{-1}}{|\Gamma|} \int_P \Phi^*g^* B \right)} 
    & & & \bigotimes_{j=p+1}^{p+q} E_{\gamma_jg} 
  }
$$
and its commutativity follows from Stokes theorem and the relation \ref{relation B},
\begin{eqnarray*}
\int_P \Phi^* g^* B - \int_P \Phi^* B & = & \int_P  \Phi^* d A_q
  =  \int_{\partial P} \Phi^* A_g \\
&  =&   \sum_{i=1}^{p} \int_{Q_i} \phi^* A_g - \sum_{j=p+1}^{p+q} \int_{Q_j} \phi^* A_g.
\end{eqnarray*} 

The previous string connection is also compatible with the connection that we have associated to the loop orbifold in \cite{LupercioUribe6}. Let's recall the construction.
The connection for us will be a linear functional $\Delta$ on the tangent
 space of the loop orbifold (a 1-form on the loop orbifold).
For $\phi: Q \to M$ with $Q$ a $\Gamma$-principal bundle and a tangent vector to it, namely 
a section $\mu : Q \to \phi^*TM$ such that $\mu$ is $\phi_\#$ equivariant we
define 
\begin{eqnarray*}
\langle \Delta_\phi, \mu \rangle :=  \frac{1}{|\Gamma|} \int_Q  \Phi^*\left( i_\mu B \right)
\end{eqnarray*}
where $i_\mu$ is contraction on the direction of $\mu$. In \cite{LupercioUribe6} we have proved that
$\Delta$ together with $F$, the gluing information of the bundle $E$ (see \ref{def. F}), form a Beilinson
Deligne cocycle over the loop orbifold, hence a line bundle with connection over it.

So we can conclude this paper with the following result
\begin{theorem}
Let $\xi$ be a global gerbe with connection over $[M/G]$ and $E$ the line
bundle with connection induced by it via transgression. Then $\xi$ permits to define
a string connection $\UU$ over the line bundle $E$ of the loop groupoid $(\LL [M/G])/G$. 
\end{theorem}

We actually have done more. We claim that we can construct a string connection
over the loop orbifold of a general orbifold (smooth Deligne-Mumford stack) from any gerbe with connection. This is the subject of the sequel to this paper.


\bibliographystyle{amsplain}
\bibliography{difchar}

\providecommand{\bysame}{\leavevmode\hbox to3em{\hrulefill}\thinspace}
\begin{thebibliography}{10}

\bibitem{AdemRuan}
A.~Adem and Y.~Ruan, \emph{Twisted orbifold {$K$}-theory}, Comm. Math. Phys.
  \textbf{237} (2003), no.~3, 533--556.

\bibitem{AndoFrench}
M.~Ando and C.~French, \emph{Discrete torsion for the supersingular orbifold
  sigma genus}, arXiv:math.AT/0308068.

\bibitem{Audin}
M.~Audin, \emph{The topology of torus actions on symplectic manifolds},
  Progress in Mathematics, vol.~93, Birkh\"auser Verlag, Basel, 1991,
  Translated from the French by the author.

\bibitem{Brylinski}
J-L. Brylinski, \emph{Loop spaces, characteristic classes and geometric
  quantization}, Progress in Mathematics, vol. 107, Birkh\"auser Boston Inc.,
  Boston, MA, 1993.

\bibitem{BrylinskiMcLaughlin}
J-L. Brylinski and D.~McLaughlin, \emph{A geometric construction of the first
  {P}ontryagin class}, Quantum topology, Ser. Knots Everything, vol.~3, World
  Sci. Publishing, River Edge, NJ, 1993, pp.~209--220.

\bibitem{Murray02}
S.~Carey, L.~Johnson and M.~Murray, \emph{Holonomy on {D}-branes},
  arXiv:hep-th/0204199.

\bibitem{CheegerSimons}
J.~Cheeger and J.~Simons, \emph{Differential characters and geometric
  invariants}, Geometry and topology (College Park, Md., 1983/84), Lecture
  Notes in Math., vol. 1167, Springer, Berlin, 1985, pp.~50--80.

\bibitem{DupontHainZucker}
J.~Dupont, R.~Hain, and S.~Zucker, \emph{Regulators and characteristic classes
  of flat bundles}, The arithmetic and geometry of algebraic cycles (Banff, AB,
  1998), CRM Proc. Lecture Notes, vol.~24, Amer. Math. Soc., Providence, RI,
  2000, pp.~47--92.

\bibitem{Esnault}
H.~Esnault, \emph{Characteristic classes of flat bundles}, Topology \textbf{27}
  (1988), no.~3, 323--352.

\bibitem{HopkinsSinger}
M.J. Hopkins and I.M. Singer, \emph{Quadratic functions in geometry,
  topology,and m-theory}, arXiv:math.AT/0211216.

\bibitem{LupercioUribe1}
E.~Lupercio and B.~Uribe, \emph{Gerbes over orbifolds and twisted k-theory},
  arXiv:math.AT/0105039, to appear in Comm. Math. Phys.

\bibitem{LupercioUribe6}
\bysame, \emph{Holonomy for gerbes over orbifolds}, arXiv:math.AT/0307114.

\bibitem{LupercioUribe3}
\bysame, \emph{Deligne cohomology for orbifolds, discrete torsion and
  b-fields}, Geometric and Topological methods for Quantum Field Theory, World
  Scientific, 2002.

\bibitem{LupercioUribe2}
\bysame, \emph{Loop groupoids, gerbes, and twisted sectors on orbifolds},
  Orbifolds in mathematics and physics (Madison, WI, 2001), Contemp. Math.,
  vol. 310, Amer. Math. Soc., Providence, RI, 2002, pp.~163--184.

\bibitem{Moerdijk98}
I.~Moerdijk, \emph{Proof of a conjecture of {A}. {H}aefliger}, Topology
  \textbf{37} (1998), no.~4, 735--741.

\bibitem{Moerdijk2002}
\bysame, \emph{Orbifolds as groupoids: an introduction}, Orbifolds in
  mathematics and physics (Madison, WI, 2001), Contemp. Math., vol. 310, Amer.
  Math. Soc., Providence, RI, 2002, pp.~205--222.

\bibitem{Ruan2}
Y.~Ruan, \emph{Discrete torsion and twisted orbifold cohomology},
  arXiv:math.AG/0005299.

\bibitem{Ruan3}
\bysame, \emph{Gerbe and twisted orbifold quantum cohomology}, preprint.

\bibitem{Ruan1}
\bysame, \emph{Stringy geometry and topology of orbifolds}, Symposium in Honor
  of C. H. Clemens (Salt Lake City, UT, 2000), Contemp. Math., vol. 312, Amer.
  Math. Soc., Providence, RI, 2002, pp.~187--233.

\bibitem{Segal02}
G.~Segal, \emph{Topological structures in string theory}, R. Soc. Lond. Philos.
  Trans. Ser. A Math. Phys. Eng. Sci. \textbf{359} (2001), no.~1784,
  1389--1398, Topological methods in the physical sciences (London, 2000).

\bibitem{Sharpe}
E.~Sharpe, \emph{Discrete torsion, quotient stacks, and string orbifolds},
  Orbifolds in mathematics and physics (Madison, WI, 2001), Contemp. Math.,
  vol. 310, Amer. Math. Soc., Providence, RI, 2002, pp.~301--331.

\bibitem{Weil52}
A.~Weil, \emph{Sur les th\'eor\`emes de de {R}ham}, Comment. Math. Helv.
  \textbf{26} (1952), 119--145.

\end{thebibliography}
\end{document}